\documentclass[11pt]{amsart}

\usepackage[a4paper,hmargin=3.5cm,vmargin=4cm]{geometry}
\usepackage{amsfonts,amssymb,amscd,amstext}
\usepackage{graphicx}
\usepackage[dvips]{epsfig}
\usepackage{quoting}
\usepackage{hyperref}

\usepackage{fancyhdr}
\input xy
\xyoption{all}

\usepackage{enumerate}
\usepackage{titlesec}
\usepackage{mathrsfs}

\titleformat{\section}
{\filcenter\bfseries\large} {\thesection{.}}{0.2cm}{}
\titleformat{\subsection}[runin]
{\bfseries} {\thesubsection{.}}{0.15cm}{}[.]
\titleformat{\subsubsection}[runin]
{\em}{\thesubsubsection{.}}{0.15cm}{}[.]

\usepackage[up,bf]{caption}

\newtheorem{theorem}{Theorem}

\newtheorem{corollary}{Corollary}

\theoremstyle{definition}
\newtheorem{definition}{Definition}

\newtheorem{problem}{Problem}

\usepackage{color}

\begin{document}


\fancyhead[LO]{The Yang problem for complete bounded complex submanifolds}
\fancyhead[RE]{A.\ Alarc\'on}
\fancyhead[RO,LE]{\thepage}

\thispagestyle{empty}



\begin{center}
{\bf\LARGE The Yang problem for complete bounded complex submanifolds: a survey
} 

\bigskip

%
%
{\large\bf Antonio Alarc\'on}
\end{center}


%
%

\bigskip

\begin{quoting}[leftmargin={5mm}]
{\small
\noindent {\bf Abstract}\hspace*{0.1cm}
We survey the history as well as recent progress in the Yang problem concerning the existence of complete bounded complex submanifolds of the complex Euclidean spaces. We also point out some open questions on the topic.


\noindent{\bf Keywords}\hspace*{0.1cm} 
Complex manifold, complete Riemannian manifold, holomorphic curve.


\noindent{\bf Mathematics Subject Classification (2020)}\hspace*{0.1cm} 
32H02, 
32B15, 
32E30. 
}
\end{quoting}


%
%
\section{The problem and the first solution}
\label{sec:problem}

In 1977, P.\ Yang asked the following (see \cite[p.\ 135, Question II]{Yang1977}):
\begin{problem}\label{qu:Yang}
Do there exist complete immersed complex submanifolds $\varphi\colon M^k\to\mathbb{C}^n$ $(1\le k<n)$ with bounded image?
\end{problem} 
Completeness is a very natural condition to impose on a Riemannian manifold when one is interested in its global properties. 
Recall that an immersed submanifold $\varphi\colon M\to\mathbb{C}^n$ is said to be complete if the Riemannian metric $\varphi^*d\sigma^2$ induced on $M$ by pulling back the Euclidean metric $d\sigma^2$ on $\mathbb{C}^n$ by the immersion $\varphi$ is a complete metric on $M$: geodesics go on indefinitely. By the Hopf-Rinow theorem, $\varphi$ is complete if and only if $\varphi\circ\gamma\colon[0,1)\to\mathbb{C}^n$ has infinite Euclidean length for every divergent path $\gamma\colon[0,1)\to M$.   Every compact Riemannian manifold is complete, but compact complex manifolds cannot be found in $\mathbb{C}^n$ by the maximum principle for holomorphic functions, so the question in Problem \ref{qu:Yang} is in order. 

The main original motivation for P.\ Yang to pose the aforementioned question is that a positive answer would prove the existence of complete immersed complex submanifolds $M^k\to\mathbb{C}^{2n}$ with strongly negative holomorphic sectional curvature \cite{Yang1977JDG}. On the other hand, since complex submanifolds of $\mathbb{C}^n$ are minimal (i.e., critical points for the volume functional; see e.g. \cite{Lawson1980,ColdingMinicozzi2006BLMS}, among many others, for an introduction to the subject), the Yang problem is also related to the so-called Calabi-Yau problem, which dates back to E.\ Calabi's conjectures from 1965 \cite[p.\ 170]{Calabi1965Conjecture} and asked whether there are complete bounded minimal hypersurfaces in $\mathbb{R}^n$ $(n\ge 3)$; we refer to \cite[Ch.\ 7]{AlarconForstnericLopez2021Book} for a recent survey on this fascinating topic. Nevertheless, the Yang problem for complex submanifolds has become an active focus of interest in its own right, having received many important contributions in the last decade; see \cite[\textsection 4.3]{Forstneric2018Survey} for a brief introduction to the topic.

We shall denote by $\mathbb{D}=\{\zeta\in\mathbb{C}\colon |\zeta|<1\}$ the open unit disc in $\mathbb{C}$ and by $\mathbb{B}_n=\{z=(z_1,\ldots,z_n)\in\mathbb{C}^n\colon |z|<1\}$ the open unit ball in $\mathbb{C}^n$ for $n\ge 2$. The first affirmative answer to the question in Problem \ref{qu:Yang} was given only two years later, in 1979, by P.W.\ Jones. Recall that a map $f\colon X\to Y$ between topological spaces is said to be proper if $f^{-1}(K)\subset X$ is compact for every compact set $K\subset Y$.
%
%
\begin{theorem}\label{th:Jones}
{\bf (Jones \cite{Jones1979PAMS})} There exist a complete bounded immersed complex disc $\mathbb{D}\to\mathbb{C}^2$, a complete bounded embedded complex disc $\mathbb{D}\hookrightarrow\mathbb{C}^3$, and a complete properly embedded complex disc $\mathbb{D}\hookrightarrow\mathbb{B}_4$. 
\end{theorem}

P.W.\ Jones' construction method is strongly complex analytic. It relies on using the BMO duality theorem in order to find a pair of bounded holomorphic functions $f_1$ and $f_2$ on $\mathbb{D}$ satisfying the property that
\[
	\int_\gamma\big( |f_1'(\zeta)|+|f_2'(\zeta)|\big)\, d\sigma(\zeta)=+\infty
\]
for all paths $\gamma\subset\mathbb{D}$ terminating on $\mathbb{S}^1=b\mathbb{D}$, where $\sigma$ denotes Euclidean arc length. It follows that $\mathbb{D}\ni\zeta\mapsto(\zeta,f_1(\zeta),f_2(\zeta))$ is a complete bounded embedded complex disc in $\mathbb{C}^3$; the other two assertions in the theorem are obtained by slight modifications of this procedure.

%
%
\section{Curves}
\label{sec:curves}

Despite having to wait more than three decades for it, Theorem \ref{th:Jones} has been generalized in several directions. The first extension of P.W.\ Jones' existence result was given by A.\ Alarc\'on and F.J. L\'opez in 2013 and concerns the topology of the examples.
\begin{theorem}\label{th:AL2013}
{\bf (Alarc\'on-L\'opez \cite{AlarconLopez2013MA})} Let $n\ge 2$. Every open orientable surface $S$ admits a complex structure $J$ such that the open Riemann surface $M=(S,J)$ carries a complete proper holomorphic immersion $M\to\mathbb{B}_n$ which is an embedding if $n\ge 3$. 
The same holds true if we replace the ball by any convex domain in $\mathbb{C}^n$.
\end{theorem}

The embeddedness condition for $n\ge 3$ in this statement was not explicitly stated in \cite{AlarconLopez2013MA}; nevertheless, it trivially follows from a standard transversality argument, as was later pointed out by A.\ Alarc\'on and F.\ Forstneri\v c in \cite{AlarconForstneric2014IM}. We emphasize that the examples in Theorem \ref{th:AL2013} may have any topological type, even infinite. 

\begin{figure}[ht]
\includegraphics[scale=.05]{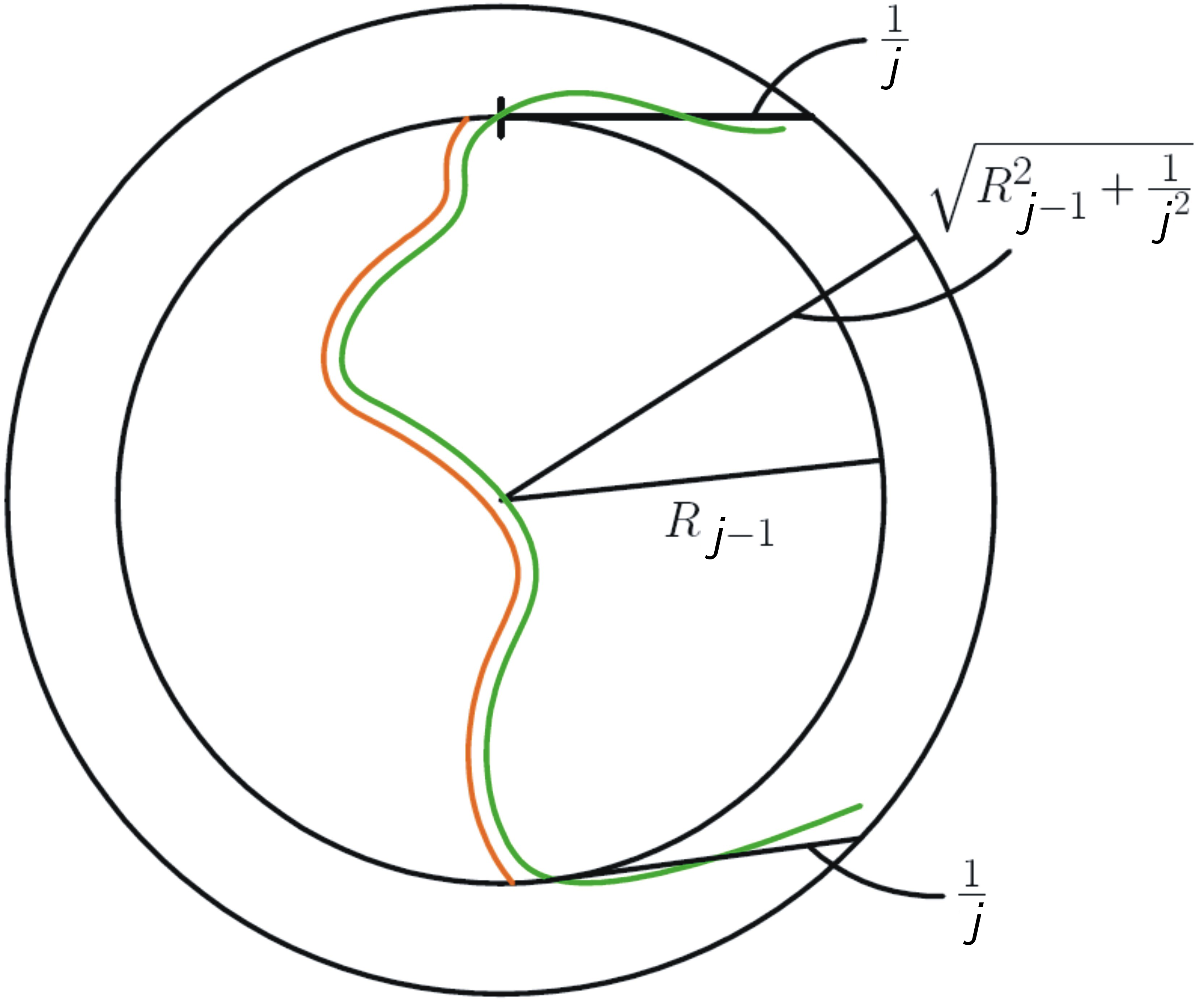}
\caption{Schematic representation of the geometry of the deformations used in the inductive construction in \cite{AlarconLopez2013MA}. The key idea is to deform the submanifold near its boundary by pushing each boundary point in a direction orthogonal to its position vector.}
\label{fig:Pi} 
\end{figure}

The proof in \cite{AlarconLopez2013MA} is completely different from that in \cite{Jones1979PAMS}. In particular, it is much more geometric, and is reminiscent of the method developed by N.\ Nadirashvili in his seminal paper \cite{Nadirashvili1996IM} for constructing a complete bounded minimal disc in $\mathbb{R}^3$; see also \cite[\textsection 7.1]{AlarconForstnericLopez2021Book}. The construction goes by induction, the rough idea being the following. In the step $j\in\mathbb{N}$ we begin with a smoothly bounded compact complex curve, say $X_{j-1}$, whose boundary $bX_{j-1}$ lies inside the ball $R_{j-1}\mathbb{B}_n$ of some radius $R_{j-1}>0$ but close to its boundary sphere. Then, we apply to $X_{j-1}$ a deformation that is arbitrarily small outside a neighborhood of $bX_{j-1}$ and pushes each boundary point $x\in bX_{j-1}$ a distance approximately $1/j$ in a direction perpendicular to the position vector of $x$ in $\mathbb{C}^n$. In this way we increase the boundary distance from a fixed interior point of the curve an amount of approximately $1/j$, while the extrinsic diameter is increased, by Pythagoras' theorem, an amount of the order of $1/j^2$. Moreover, we ensure that the boundary of the new complex curve, $X_j$, lies inside the ball of radius 
\[
	R_j=\sqrt{R_{j-1}^2+1/j^2}>R_{j-1}
\]
but close to its boundary sphere. See Figure \ref{fig:Pi}. 
%
%
Since 
\[
	\sum_{j\ge 1}1/j=\infty\quad \text{and}\quad \sum_{j\ge 1}1/j^2<\infty,
\] 
if we arrange this process in the right way then we obtain in the limit a complete complex curve that is contained in the ball of radius $R=\lim_{j\to\infty}R_j<\infty$ and is proper in it. In order to prescribe the topology of the curve, we begin with a disc and at each step of the inductive construction apply a surgery which enables us to add  either a handle or a boundary component to a given compact bordered complex curve in $\mathbb{C}^n$. Finally, if $n\ge 3$, then a general position argument allows to guarantee that all complex curves in the sequence are embedded, and hence the limit one can be granted to be embedded as well. The main tool in order to make the described deformations is the theory of uniform approximation for holomorphic functions on open Riemann surfaces; in particular, the Runge-Mergelyan theorem (see E. Bishop \cite{Bishop1958PJM} or e.g. Theorem 5 in the survey on holomorphic approximation by J.E.\ Forn{\ae}ss, F.\ Forstneri\v{c}, and E.F.\ Wold \cite{FornaesForstnericWold2020}).

A different construction of complete bounded complex curves in $\mathbb{C}^2$ can be found in \cite{MUY2}, where F. Mart\'in, M. Umehara, and K. Yamada gave examples with arbitrary finite genus and finitely many ends. Their technique relies on the existence of a simply connected complete bounded holomorphic null curve in $\mathbb{C}^3$ (such a curve was first constructed in \cite{AlarconLopez2013MA}, an alternative construction was given later by L.\ Ferrer, F.\ Mart\'in, M.\ Umehara, and K.\ Yamada in \cite{FerrerMartinUmeharaYamada2014}; see the mentioned sources or e.g. \cite[\textsection 2.3]{AlarconForstnericLopez2021Book} for definitions) and modifies a method developed by F.J.\ L\'opez in \cite{Lopez1998TAMS} for constructing complete minimal surfaces in $\mathbb{R}^3$ of hyperbolic conformal type.

Despite the flexibility of the methods in \cite{AlarconLopez2013MA,MUY2}, none of them allows to control the complex structure on the curve, except of course in the simply connected case: every bounded immersed complex curve in $\mathbb{C}^n$ must have hyperbolic conformal type, and hence if it is simply connected it is biholomorphic to the disc $\mathbb{D}$. Indeed, since the construction in \cite{AlarconLopez2013MA} relies on Runge's theorem, at a certain stage one does not have enough information on the placement in $\mathbb{C}^n$ of some parts of the curve, and hence one is forced to cut away some small pieces of the curve to ensure its boundedness. This makes one to lose the control of the complex structure of the curve. This difficulty was overcome by A.\ Alarc\'on and F.\ Forstneri\v c in \cite{AlarconForstneric2013MA}, also published in 2013, where two additional complex analytic tools were introduced into the game, namely, the method of F.\ Forstneri\v c and E.F.\ Wold \cite{ForstnericWold2009} for exposing boundary points of a complex curve in $\mathbb{C}^n$ and the use of approximate solutions to Riemann-Hilbert boundary value problems. (The former is a modern technique that has led to important progress in the classical Forster-Bell-Narasimhan Conjecture asking whether every open Riemann surface admits a proper holomorphic embedding into $\mathbb{C}^2$ \cite{BellNarasimhan1990EMS,Forster1970CMH}; we refer to F.\ Forstneri\v c \cite[\textsection 9.10-9.11]{Forstneric2017E} for a survey on this long-standing open problem. On the other hand, the use of the Riemann-Hilbert problem for constructing proper holomorphic maps has a long history; we refer to F.\ Forstneri\v c and J.\ Globevnik \cite{ForstnericGlobevnik1992CMH}, B.\ Drinovec Drnov\v sek and F.\ Forstneri\v c \cite{DrinovecForstneric2007DMJ}, and the references therein.) The implementation of these new tools enabled to substantially  simplify the construction in \cite{AlarconLopez2013MA} and, moreover, to control the complex structure on the curve. Recall that a bordered Riemann surface is an open connected Riemann surface $M$ that is the interior, $M=\overline M\setminus b\overline M$, of a compact one dimensional complex manifold $\overline M$ with smooth boundary $b\overline M$ consisting of finitely many closed Jordan curves; such an $\overline M$ is called a compact bordered Riemann surface.
%
%
\begin{theorem}\label{th:AF2013}
{\bf (Alarc\'on-Forstneri\v c \cite{AlarconForstneric2013MA})} Let $n\ge 2$ be an integer. Every bordered Riemann surface $M$ admits a complete proper holomorphic immersion $M\to\mathbb{B}_n$ that is an embedding if $n\ge 3$.
The same holds true if we replace the ball by any pseudoconvex domain in $\mathbb{C}^n$.
\end{theorem}

More generally, it is shown in \cite{AlarconForstneric2013MA} that if $X$ is a Stein manifold of dimension $n\ge 2$ endowed with a hermitian metric, then every bordered Riemann surface $M$ admits a complete proper holomorphic immersion $M\to X$ that can be chosen an embedding if $n\ge 3$. Recall that a Stein manifold is the same thing as a closed complex submanifold of a complex Euclidean space; we refer to \cite{Forstneric2017E} for a comprehensive monograph on the theory of Stein manifolds. A domain $D$ in $\mathbb{C}^n$ $(n\ge 2)$ is a Stein manifold if and only if it is pseudoconvex, meaning that there is a strictly plurisubharmonic exhaustion function $D\to\mathbb{R}$. This happens if and only if $D$ is a domain of holomorphy and if and only if $D$ is holomorphically convex. For instance, every convex domain in $\mathbb{C}^n$ is pseudoconvex. We refer to the monographs by R.M.\ Range \cite{Range1986} (see also the introductory note \cite{Range2012}) and L.\ H{\"o}rmander \cite{Hormander1990,Hormander1994}  for background on the subject.

More recently, the construction technique in \cite{AlarconForstneric2013MA} has been refined to produce complete bounded complex curves with control on the complex structure and with some further control on the asymptotic behavior. In particular, there are such curves bounded by Jordan curves. The following is a compilation of results by A.\ Alarc\'on, I.\ Castro Infantes, B. Drinovec Drnov\v sek, F.\ Forstneri\v c, and F.J. L\'opez \cite{AlarconLopez2013IJM,AlarconDrinovecForstnericLopez2015PLMS,AlarconCastroInfantes2018GT,AlarconForstneric2021RMI,Forstneric2022RMI}.
\begin{theorem}\label{th:Jordan}
Let $R$ be a compact Riemann surface and assume that $M=R\setminus\bigcup_{i\in I} D_i$ is a domain in $R$ whose complement is a finite or countable union of pairwise disjoint, smoothly bounded closed discs (i.e., diffeomorphic images of $\overline{\mathbb{D}}$). The following assertions hold true: 
\begin{enumerate}[1.]
\item $M$ is the complex structure of a complete bounded complex curve in $\mathbb{C}^2$. In fact, for any $n\ge 2$ there is a continuous map $\varphi\colon \overline M\to\mathbb{C}^n$ such that the restricted map $\varphi|_M\colon M\to\mathbb{C}^n$ is a complete holomorphic immersion (embedding if $n\ge 3$) and $\varphi|_{bM}\colon bM=\bigcup_{i\in I} bD_i\to\mathbb{C}^n$ is injective. {\bf (Alarc\'on-Drinovec~Drnov{\v{s}}ek-Forstneri\v c-L\'opez \cite{AlarconDrinovecForstnericLopez2015PLMS}, Alarc\'on-Forstneri\v c \cite{AlarconForstneric2021RMI}.)}
\smallskip
\item If $I$ is finite, then there is a continuous map $\varphi\colon \overline{M}\to\overline{\mathbb{B}}_n$ $(n\ge 3)$ such that $\varphi(bM)\subset b\mathbb{B}_n$ and $\varphi|_M\colon M\to\mathbb{B}_n$ is a complete proper holomorphic immersion (embedding if $n\ge 3$). The same holds true with $\mathbb{B}_n$ replaced by any convex domain in $\mathbb{C}^n$. {\bf (Alarc\'on-Drinovec~Drnov{\v{s}}ek-Forstneri\v c-L\'opez \cite{AlarconDrinovecForstnericLopez2015PLMS}.)}
\smallskip 
\item If I is finite, then for any domain $D$ in $\mathbb{C}^n$ there is a complete holomorphic immersion $\varphi\colon M\to D$ whose image is a dense subset of $D$. If $n\ge 3$ then $\varphi$ can be chosen injective. {\bf (Alarc\'on-Castro Infantes \cite{AlarconCastroInfantes2018GT}.)}
\smallskip 
\item There is a Cantor set $C$ in $R$ whose complement admits a complete holomorphic immersion $R\setminus C\to\mathbb{C}^2$ with bounded image. There also exist a Cantor set $C$ in $R$ and a complete holomorphic embedding $R\setminus C\hookrightarrow\mathbb{C}^3$ with bounded image. {\bf (Forstneri\v c \cite{Forstneric2022RMI}.)}
\smallskip
\item There is a Cantor set $C$ in $M$ and a continuous map $\varphi\colon \overline{M}\setminus C\to \mathbb{C}^n$ $(n\ge 2)$ such that $\varphi|_{M\setminus C}\colon M\setminus C\to\mathbb{C}^n$ is a complete holomorphic immersion and $\varphi|_{bM}\colon bM=\bigcup_{i\in I} bD_i\to\mathbb{C}^n$ is injective. If $n\ge 3$ then $C$ can be chosen so that $\varphi\colon \overline{M}\setminus C\to \mathbb{C}^n$ is an injective map. {\bf (Forstneri\v c \cite{Forstneric2022RMI}.)}
\end{enumerate}
\end{theorem}

Summarizing, by the year 2015 there were available in the literature several constructions of complete bounded complex curves immersed  in $\mathbb{C}^2$ and embedded in $\mathbb{C}^3$, allowing a high control on the asymptotic behavior (proper in the ball or in a given pseudoconvex domain, bounded by Jordan curves, etc), on the topology, and on the complex structure of the examples. However, the construction of complete bounded embedded complex curves in $\mathbb{C}^2$ turns out to be a much more challenging undertaking, and the question whether such curves exist remained open (see \cite[Question 1]{AlarconForstneric2013MA}). Recall that complex curves are generically embedded in $\mathbb{C}^n$ for $n\ge 3$, meaning that self-intersections can be removed by applying small deformations, while self-intersections of complex curves in $\mathbb{C}^2$, which are generically double points, are stable under such deformations. That is the main reason why the task is a more difficult one.

%
%

The Yang problem for embedded complex curves in the affine plane $\mathbb{C}^2$ was finally settled by A.\ Alarc\'on and F.J.\ L\'opez in a paper published in 2016.
%
%
\begin{theorem}\label{th:AL2016}
{\bf (Alarc\'on-L\'opez \cite{AlarconLopez2016JEMS})} 
Every convex domain in $\mathbb{C}^2$ admits complete properly embedded complex curves.
\end{theorem}

The proof goes by induction and involves an approximation process by embedded complex curves in $\mathbb{C}^2$. The main step in the construction is to prove that every compact embedded complex curve $X$ in $\mathbb{C}^2$ with the boundary $bX$ lying in the boundary $bD$ of a regular strictly convex domain $D\subset\mathbb{C}^2$ may be approximated by another compact embedded complex curve $\tilde X$ with $b\tilde X\subset b\tilde D$, for any given convex domain $\tilde D\subset\mathbb{C}^2$ with $\overline D\subset\tilde D$. The new curve $\tilde X$ is ensured to contain a biholomorphic copy of $X$, which we denote by $X$ as well, and the main point is to guarantee that $\tilde X\setminus X\subset\tilde D\setminus D$ and the intrinsic Euclidean distance in $\tilde X$ from $X$ to $b\tilde X$ is suitably larger (in a Pythagorical way similar to that explained in Figure \ref{fig:Pi}) than the distance from $D$ to $b\tilde D$. These conditions are the key for obtaining embeddedness, completeness, and properness of the limit complex curve in the limit convex domain. In order to guarantee the embeddedness of $\tilde X$ a standard self-intersection removal method consisting of replacing every normal crossing in an immersed complex curve in $\mathbb{C}^2$ by an embedded annulus is applied. This surgery may generate shortcuts in the arising desingularized curve $\tilde X$, thereby giving rise to divergent paths of shorter length. This is an important difficulty for ensuring completeness; for instance, if one applies this surgery at each step in the inductive construction in the proof of Theorem \ref{th:AL2013} or Theorem \ref{th:AF2013}, then one still obtains in the limit a properly embedded complex curve in $\mathbb{B}_2$ (or, more generally, in any given pseudoconvex domain of $\mathbb{C}^2$), but it need not be complete. A main novelty in the construction in \cite{AlarconLopez2016JEMS} is provided a good enough (say, in a Pythagorical sense) estimate of the growth of the intrinsic Euclidean diameter of the desingularized, embedded complex curve $\tilde X$; this is achieved by keeping a stronger control on the placement of $\tilde X\setminus X$ in $\tilde D\setminus D$. The main tool in this construction continues to be the classical Runge-Mergelyan approximation theorem for holomorphic functions on open Riemann surfaces.

The aforementioned surgery may increase the topological genus of the curve, and so there is no control on the topology of the examples in Theorem \ref {th:AL2016}. They actually seem to have infinite genus, and can be ensured to have infinite topology. It therefore remained an open question whether there are complete bounded embedded complex curves in $\mathbb{C}^2$ of finite topology (see  \cite[Question 1.5]{AlarconLopez2016JEMS}).

%
%
\section{Submanifolds of arbitrary dimension}
\label{sec:higher}

All examples of complete bounded complex submanifolds we have discussed so far are of complex dimension one; i.e., complex curves in $\mathbb{C}^n$. The first known such submanifolds of higher dimension were also given in \cite{AlarconForstneric2013MA}.
%
%
\begin{corollary}\label{co:AF2013}
{\bf (Alarc\'on-Forstneri\v c \cite{AlarconForstneric2013MA})} If $D$ is a relatively compact, strongly pseudoconvex domain in a Stein manifold of dimension $k\ge 1$, then there are a complete proper holomorphic immersion $D\to (\mathbb{B}_2)^{2k}\subset\mathbb{C}^{4k}$ and a complete proper holomorphic embedding $D\hookrightarrow (\mathbb{B}_3)^{2k+1}\subset\mathbb{C}^{6k+3}$.
\end{corollary}

The corollary is obtained by the following simple trick which was pointed out by J.E.\ Forn\ae ss. In \cite{DrinovecForstneric2010AJM}, B.\ Drinovec Drnov\v sek and F.\ Forstneri\v c proved that every domain $D$ as in the statement admits a proper holomorphic immersion $g\colon D\to \mathbb{D}^{2k}$ into the polydisc $\mathbb{D}^{2k}=\mathbb{D}\times \stackrel{2k}{\ldots}\times\mathbb{D}\subset\mathbb{C}^{2k}$. Choose a complete proper holomorphic immersion $\varphi\colon \mathbb{D}\to\mathbb{B}_2$ provided by Theorem \ref{th:AF2013}, and consider the proper holomorphic map $\varphi^{2k}\colon \mathbb{D}^{2k}\to (\mathbb{B}_2)^{2k}$ given by
\[
	\varphi^{2k}(\zeta_1,\ldots,\zeta_{2k})=(\varphi(\zeta_1),\ldots,\varphi(\zeta_{2k})),\quad (\zeta_1,\ldots,\zeta_{2k})\in  \mathbb{D}^{2k}.
\]
It is then easily checked that $\varphi^{2k}\circ g\colon D\to (\mathbb{B}_2)^{2k}$ is a complete proper holomorphic immersion. A slight modification of this argument using a proper holomorphic embedding $D\hookrightarrow \mathbb{D}^{2k+1}$ (existence of such is also proved in \cite{DrinovecForstneric2010AJM}) provides a proper holomorphic embedding $D\hookrightarrow (\mathbb{B}_3)^{2k+1}$. This same trick together with Theorem \ref{th:AL2016} allows to prove the following.
%
%
\begin{corollary}\label{co:AL2016}
{\bf (Alarc\'on-L\'opez \cite{AlarconLopez2016JEMS})} 
For any $k\in\mathbb{N}$ there is a complete bounded embedded $k$-dimensional complex submanifold $M^k\hookrightarrow \mathbb{C}^{2k}$.
\end{corollary}

The ad hoc construction of the high dimensional examples in Corollaries \ref{co:AF2013} and \ref{co:AL2016} seemed to give very particular solutions to the Yang problem for complex submanifolds of dimension $\ge 2$. This led to some new questions, as whether the dimension $2k$ in Corollary \ref{co:AL2016} is optimal, and exposed the need of looking for new construction methods other than those based on the existence of complete bounded complex curves. For instance, the natural question whether the ball $\mathbb{B}_k$ $(k\ge 2)$ admits a complete proper holomorphic embedding into the ball $\mathbb{B}_n$ for some $n>k$, and, in particular, for $n=k+1$ appeared; see \cite[Question 3]{AlarconForstneric2013MA} and \cite[Question 13.2]{Globevnik2015AM}. An affirmative answer to this question in the case of sufficiently high codimension was given in 2015 by B.~Drinovec~Drnov{\v{s}}ek.
%
%
\begin{theorem}\label{th:Barbara}
{\bf (Drinovec~Drnov{\v{s}}ek \cite{Drinovec2015JMAA})}
Every bounded strictly pseudoconvex domain $D\subset\mathbb{C}^k$ $(k\in\mathbb{N})$ with $\mathcal{C}^2$ boundary admits a complete proper holomorphic embedding $D\hookrightarrow\mathbb{B}_n$ for any large enough $n\in\mathbb{N}$.
\end{theorem}

The proof in \cite{Drinovec2015JMAA} continues to use the geometric idea of deforming a compact submanifold near the boundary in orthogonal directions to the position vector (see Figure \ref{fig:Pi}), but it exploits different tools. The main new ingredients are holomorphic peak functions, going back to ideas of M.\ Hakim and N.\ Sibony \cite{HakimSibony1982} and E.\ L\o w \cite{Low1982}, and the construction of inner functions on the ball, as well as J.E.\ Forn\ae ss embedding theorem \cite{Fornaess1976AJM}. The construction relies on suitably modifying earlier methods by F.\ Forstneri\v c \cite{Forstneric1986TAMS} and E.\ L\o w \cite{Low1985} for constructing proper holomorphic maps from a strictly convex domain with $\mathcal{C}^2$ boundary in $\mathbb{C}^k$ into a unit ball of some Euclidean space of higher dimension, in order to make them complete. The construction  method in \cite{Drinovec2015JMAA} requires a sufficiently high codimension, and hence the following question remains open; see \cite[Question 3]{AlarconForstneric2013MA} and \cite[Question 13.2]{Globevnik2015AM}.
\begin{problem}
Does there exist a complete proper holomorphic embedding $\mathbb{B}_k\hookrightarrow \mathbb{B}_{k+1}$ for $k\ge 2$?
\end{problem}

%
%
\section{Hypersurfaces}
\label{sec:hypersurfaces}

%
%

Recall that complex submanifolds of dimension $k\in\mathbb{N}$ in $\mathbb{C}^n$ are generically embedded for $n\ge 2k+1$ (meaning that one can get rid of their self-intersections by applying small deformations on compact pieces), while for $n\le 2k$ self-intersections of $k$-dimensional complex submanifolds  in $\mathbb{C}^n$ are stable under small deformations. This is the main reason why the Yang problem for embedded submanifolds is much more difficult in low codimension; in particular, for hypersurfaces. In the lowest dimensional case of $n=2$ (i.e., for complex curves in $\mathbb{C}^2$), this question was first solved by A.\ Alarc\'on and F.J.\ L\'opez in \cite{AlarconLopez2016JEMS} (see Theorem \ref{th:AL2016} above), but the question whether there are complete bounded embedded complex hypersurfaces in $\mathbb{C}^n$, or even whether there are such complex submanifolds of dimension $k$ with $2k\ge n$, remained open for every $n\ge 3$. Indeed, note that the codimension in all examples of complete bounded embedded complex submanifold of dimension $\ge 2$ which have been mentioned so far is high. 

It was J.\ Globevnik who, in a pair of landmark papers in 2015-2016, positively settled Yang's question in Problem \ref{qu:Yang} for embeddings in arbitrary dimension and codimension; in particular, for hypersurfaces. 
%
%
\begin{theorem}\label{th:PikoAM}
{\bf (Globevnik \cite{Globevnik2015AM,Globevnik2016MA})} 
For any pair of integers $1\le k<n$ there is a complete closed embedded $k$-dimensional complex submanifold of $\mathbb{B}_n$. In particular, $\mathbb{B}_n$ admits a complete properly embedded complex hypersurface. 

The same holds true if we replace the ball by any pseudoconvex domain in $\mathbb{C}^n$.
\end{theorem}

Globevnik's approach is completely different from any previous method used in the study of the Yang problem, and it was a major breakthrough in this topic. In particular, his construction of a complete closed complex hypersurface in a given pseudoconvex domain $D\subset\mathbb{C}^n$ is implicit: the examples are obtained as level sets of highly oscillating holomorphic functions on $D$. (The existence of complete closed complex submanifolds of any higher codimension is then an obvious consequence.) To be more precise, Globevnik proved the following.
%
%
\begin{theorem}\label{th:PikoMA}
{\bf (Globevnik \cite{Globevnik2015AM,Globevnik2016MA})} 
For any pseudoconvex domain $D\subset\mathbb{C}^n$ $(n\ge 2)$ there is a holomorphic function on $D$ whose real part is unbounded above on every divergent path $\gamma\colon [0,1)\to D$ of finite length.
\end{theorem}

If a function $f\colon D\to\mathbb{C}$ on a domain $D\subset\mathbb{C}^n$ is holomorphic and nonconstant, then all its nonempty level sets $f^{-1}(c)=\{z\in D\colon f(z)=c\}$ $(c\in\mathbb{C})$ are closed complex hypersurfaces of $D$, possibly with singularities. Nevertheless, by Sard's theorem most of them are smooth, and hence properly embedded. So, if $c\in\mathbb{C}$ is such that $f^{-1}(c)\neq\varnothing$, then every divergent path $\gamma\colon [0,1)\to f^{-1}(c)$ diverges on $D$ as well. If $D$ is pseudoconvex and the function $f$ is as in Theorem \ref{th:PikoMA} then, since $\Re(f)$ is constant (and hence bounded) on $\gamma$, we have that $\gamma$ has infinite length, and $f^{-1}(c)$ is thus complete. Therefore, the level sets $f^{-1}(c)$ $(c\in\mathbb{C})$ of $f$ form a (possibly singular) holomorphic foliation of $D$ by complete closed complex hypersurfaces. Theorem \ref{th:PikoMA} thus implies the following corollary which, in turn, implies Theorem \ref{th:PikoAM}.
%
%
\begin{corollary}\label{co:PikoMA}
{\bf (Globevnik \cite{Globevnik2015AM,Globevnik2016MA})} 
Every pseudoconvex domain $D\subset\mathbb{C}^n$ $(n\ge 2)$ admits a (possibly singular) holomorphic foliation by complete closed complex hypersurfaces (most of which are smooth). 
\end{corollary}

\begin{figure}[ht]
\includegraphics[scale=.2]{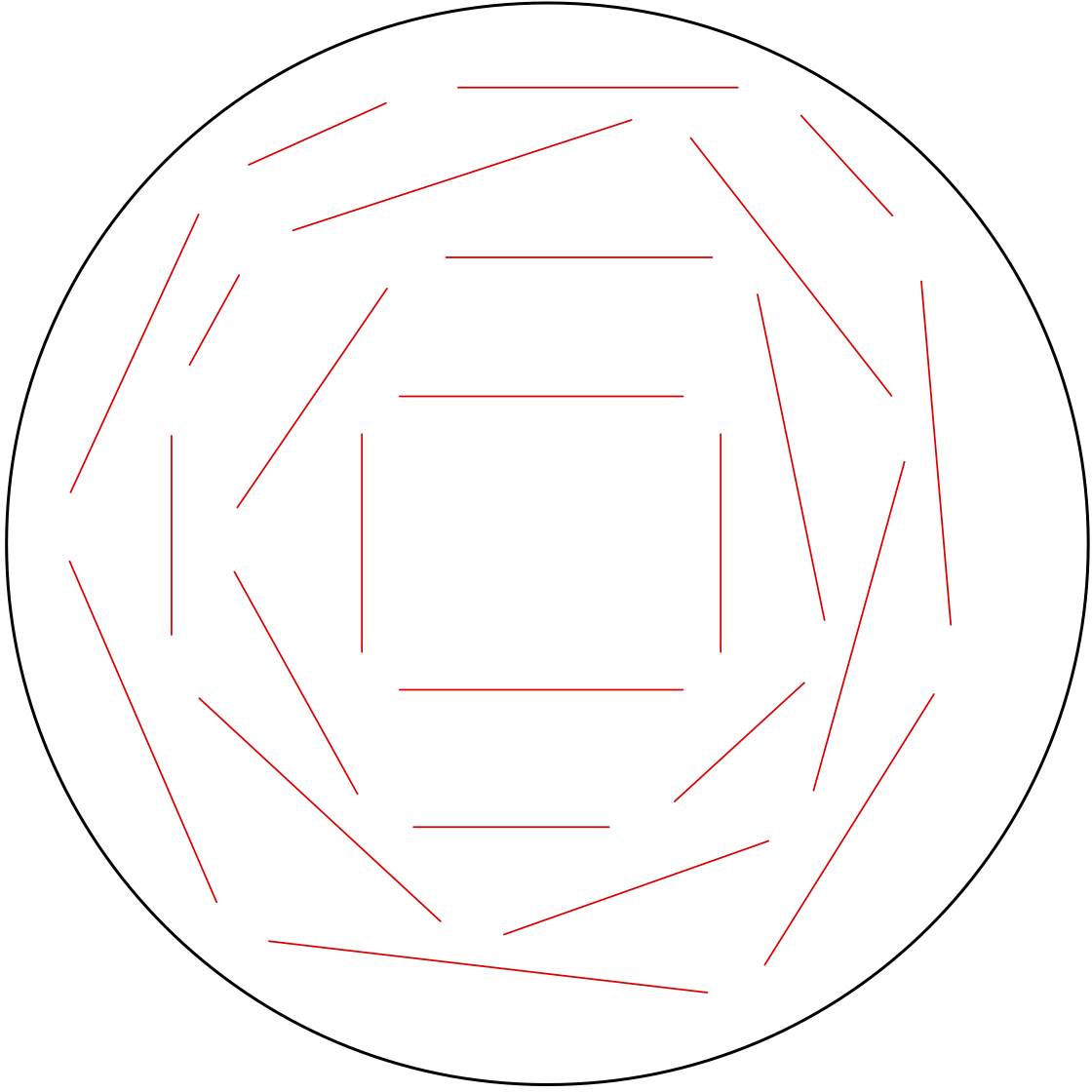}
\caption{Globevnik's labyrinth. Each piece of the labyrinth is a compact convex set in a real affine hyperplane. Every divergent path in the ball meeting at most finitely many components of the labyrinth has infinite length.}
\label{fig:L-G} 
\end{figure}

Let us outline the proof of Theorem \ref{th:PikoMA} in the case of $D=\mathbb{B}_n$ given in \cite{Globevnik2015AM}. Recall that a convex polytope $P$ in $\mathbb{R}^d$, $d\ge 2$, is a compact convex set which is the intersection of finitely many closed half-spaces. A face of $P$ is a closed convex subset $F\subset P$ such that every closed segment in $P$ whose relative interior intersects $F$ is contained in $F$. The boundary $bP$ of $P$ is the union of its faces of dimension $d-1$, while the skeleton ${\rm skel}(P)$ of $P$ is the union of all $(d-2)$-dimensional faces of $P$. Most of the work in the proof consists of constructing a sequence of convex polytopes $P_j$ in $\mathbb{C}^n=\mathbb{R}^{2n}$ and positive numbers $\theta_j$ $(j\in\mathbb{N})$ satisfying the following conditions:
\begin{enumerate}[(ii)]
\item[(i)] $P_1\subset {\rm Int}(P_2)\subset P_2\subset {\rm Int}(P_3)\subset\cdots\subset\bigcup_{j\in\mathbb{N}}P_j=\mathbb{B}_n$.
\smallskip
\item[(ii)] Denote by $U_j$ the $\theta_j$-neighborhood of ${\rm skel}(P_j)$ in $bP_j$, and set $V_j=(bP_j)\setminus U_j$, $j\in\mathbb{N}$. If $\gamma\colon[0,1)\to \mathbb{B}_n$ is a divergent path such that $\gamma([0,1))\cap V_j=\varnothing$ for all $j\ge j_0$ for some $j_0\in\mathbb{N}$, then $\gamma$ has infinite length.
\end{enumerate}
Each set $V_j$ $(j\in\mathbb{N})$ is compact and its connected components are closed convex sets in hyperplanes of $\mathbb{R}^{2n}$. The union 
\begin{equation}\label{eq:Lab}
	L=\bigcup_{j\in\mathbb{N}} V_j 
\end{equation}
of all of them is a sort of labyrinth of compact connected $(2n-1)$-dimensional convex sets in $\mathbb{B}_n\subset\mathbb{R}^{2n}$ with the property that every divergent  path $\gamma\colon[0,1)\to\mathbb{B}_n$ meeting at most finitely many components of $L$ has infinite length. See Figure \ref{fig:L-G}. 
%
%
The construction of the labyrinth in \cite{Globevnik2015AM} is very involved and belongs to convex geometry. With the labyrinth in hand, to complete the proof of Theorem \ref{th:PikoMA} an idea of J.\ Globevnik and E.L.\ Stout from \cite{GlobevnikStout1982JAM} is used in order to construct, via Runge's theorem, a sequence of holomorphic polynomials $f_j\colon \mathbb{C}^n\to\mathbb{C}$, $j\in\mathbb{N}$, such that the following conditions hold for each $j\in\mathbb{N}$:
\begin{enumerate}[(iii)]
\item[(iii)] $\Re(f_j(z))\ge j+1$ for all $z\in V_j$.
\smallskip
\item[(iv)] $|f_{j+1}(z)-f_j(z)|\le 1/2^{j+1}$ for all $z\in P_j$.
\end{enumerate}
Conditions (i) and (iv) ensure that there is a limit holomorphic function 
\[
	f=\lim_{j\to\infty}f_j\colon\mathbb{B}_n\to\mathbb{C}. 
\]
By conditions (iii) and (iv), we have $|f(z)-f_j(z)|<1$ and hence $\Re(f(z))\ge j$ for all $z\in V_j$, $j\in\mathbb{N}$. Since every divergent path $\gamma\colon[0,1)\to\mathbb{B}_n$ of finite length meets $V_j$ for infinitely many $j$'s by condition (ii), we obtain that $\Re(f)$ is unbounded above on every such $\gamma$, proving Theorem \ref{th:PikoMA} in the case of $D=\mathbb{B}_n$. 

The proof of the theorem for an arbitrary pseudoconvex domain in $\mathbb{C}^n$ $(n\ge 2)$ given in \cite{Globevnik2016MA} relies on a modification of this technique, exploiting the well known fact that every such domain is biholomorphic to a properly embedded complex submanifold of a complex Euclidean space; see \cite[Theorem 5.3.9]{Hormander1990}. (In fact, if $n\ge 2$ then every $n$-dimensional Stein manifold admits a proper holomorphic embedding into $\mathbb{C}^N$ for $N=[\frac{3n}{2}]+1$; see Y.~Eliashberg and M.~Gromov \cite{EliashbergGromov1992AM} for even $n$ and J.~Sch{\"u}rmann \cite{Schurmann1997MA} for odd $n\ge 3$. The dimension $[\frac{3n}{2}]+1$ in this assertion is the lowest prossible for all $n\ge 2$ due to purely topological reasons, as O.\ Forster pointed out by examples in \cite{Forster1970CMH}.)

Theorem \ref{th:PikoMA} also holds for $n=1$ the proof being much easier. In fact, it is
not hard to show, as an application of the classical Runge theorem in one variable, that every domain $\Omega\subset\mathbb{C}$ admits a holomorphic function whose real part is unbounded both above and below on every divergent path $[0, 1) \to\Omega$ of finite length (recall that every domain in $\mathbb{C}$ is pseudoconvex).  In this direction, P.W.\ Jones proved back in 1980 that the unit disc $\mathbb{D}$ admits a bounded holomorphic function $h$  such that all level sets of $|h|\colon \mathbb{D}\to\mathbb[0,\infty)$ have infinite length \cite{Jones1980MMJ}. Obviously, the functions in Theorem \ref{th:PikoMA} are not bounded, hence the following remains an open (but likely difficult) question.
%
%
\begin{problem}
Does the open unit ball $\mathbb{B}_n$ $(n\ge 2)$ admit a bounded holomorphic function whose level sets are all complete?
\end{problem}
%
%
%

%
%
\subsection{Controlling the topology of the examples}\label{sec:topology} 

Neither the construction in \cite{AlarconLopez2016JEMS} of complete properly embedded complex curves in the ball $\mathbb{B}_2$ nor that in \cite{Globevnik2015AM,Globevnik2016MA} of complete properly embedded complex hypersurfaces in $\mathbb{B}_n$ for arbitrary $n\ge 2$ provide any control on the topology of the examples, which, in principle, could be very complicated. In particular, it remained an open question whether there are complete properly embedded complex curves of finite topology in the ball of $\mathbb{C}^2$, and, more precisely, whether there are simply connected ones; see \cite[Question 1.5]{AlarconLopez2016JEMS} and \cite[Question 13.1]{Globevnik2015AM}. 

The first construction of complete properly embedded complex hypersurfaces in the ball with control of the topology was carried out shortly after by A.\ Alarc\'on, J.\ Globevnik, and F.J.\ L\'opez, who proved the following.
%
%
\begin{theorem}\label{th:AGL}
{\bf (Alarc\'on-Globevnik-L\'opez \cite{AlarconGlobevnikLopez2019Crelle})} Assume that $Z\subset\mathbb{C}^n$ $(n\ge 2)$ is a smooth closed complex hypersurface with  $Z\cap\mathbb{B}_n\neq\varnothing$ and let $K\subset Z\cap\mathbb{B}_n$ be a connected compact set. Then there exist a pseudoconvex Runge domain $X\subset Z$ such that $K\subset X$ and a complete proper holomorphic embedding $X\hookrightarrow \mathbb{B}_n$.

In particular, $\mathbb{B}_n$ admits smooth complete closed complex hypersurfaces which are biholomorphic to a pseudoconvex Runge domain in $\mathbb{C}^{n-1}$.
\end{theorem}

In the case of $n=2$, Theorem \ref{th:AGL} trivially implies the following more precise result. Recall that all open orientable surfaces can be realized as Runge domains of properly embedded complex curves in $\mathbb{C}^2$; see M.\ {\v{C}}erne and F.\ Forstneri\v c \cite{CerneForstneric2002} for the finite topological case and A.\ Alarc\'on and F.J.\ L\'opez \cite{AlarconLopez2013JGEA} for the arbitrary one.
%
%
\begin{corollary}\label{co:AGL}
{\bf (Alarc\'on-Globevnik-L\'opez \cite{AlarconGlobevnikLopez2019Crelle})} The open unit ball $\mathbb{B}_2$ of $\mathbb{C}^2$ admits complete properly embedded complex curves with any finite topology. 

In particular, there is a complete proper holomorphic embedding $\mathbb{D}\hookrightarrow\mathbb{B}_2$.
\end{corollary}

This settled in the positive the aforementioned questions \cite[Question 1.5]{AlarconLopez2016JEMS} and \cite[Question 13.1]{Globevnik2015AM}. The construction method in \cite{AlarconGlobevnikLopez2019Crelle} is different from the ones in \cite{AlarconLopez2016JEMS,Globevnik2015AM,Globevnik2016MA}. In particular, it introduces a new tool into the game: the use of holomorphic automorphisms of $\mathbb{C}^n$. In order to prove Theorem \ref{th:AGL} one begins with the given properly embedded hypersurface $Z$ in $\mathbb{C}^n$ and the natural embedding $Z\hookrightarrow \mathbb{C}^n$ given by the inclusion map. Then, in a recursive way, one composes this initial embedding with a sequence of holomorphic automorphisms $\Phi_j$ of $\mathbb{C}^n$, $j\in\mathbb{N}$, which converges uniformly on compact subsets of $\mathbb{B}_n$. One obtains in this way a sequence of proper holomorphic embeddings $\Psi_j\colon Z\hookrightarrow\mathbb{C}^n$, where $\Psi_j=\Phi_j\circ\Phi_{j-1}\circ\cdots\circ\Phi_1$, $j\in\mathbb{N}$, whose images converge uniformly on compact subsets of $\mathbb{B}_n$ to a properly embedded complex hypersurface of $\mathbb{B}_n$. If this process is carried out in the right way, then one can guarantee that a connected component $V$ of the limit hypersurface is biholomorphic to a pseudoconvex Runge domain $X$ in $Z$ which is closed in $\mathbb{B}_n$ and contains $K$. The main point is that the sequence of holomorphic automorphisms $\Phi_j$ of $\mathbb{C}^n$ can be chosen to make sure that $V$ is disjoint from a labyrinth $L$ of compact sets in $\mathbb{B}_n$ as that constructed by J.\ Globevnik in \cite{Globevnik2015AM} (see \eqref{eq:Lab} and Figure \ref{fig:L-G}), and hence $V$ is in addition complete. Indeed, denote by $\{L_j\}_{j\in\mathbb{N}}$ the family of connected components of $L$ and note that, up to reordering, the following properties are satisfied for all $j\in\mathbb{N}$:
\begin{enumerate}[(b)]
\item[(a)] $L_j$ is a compact convex set contained in an affine real hyperplane $H_j$ of $\mathbb{C}^n$.
\smallskip
\item[(b)] There exists a compact convex set $K_j\subset\mathbb{C}^n$ such that $\bigcup_{i=1}^{j-1} L_i\subset K_j$ and $K_j\cap H_l=\varnothing$ for all $l\ge j$.
\end{enumerate}
Moreover, up to removing finitely many pieces of the labyrinth if necessary, we may also assume that:
\begin{enumerate}[(a)]
\item[(c)] There is a compact convex set $K_0\subset\mathbb{C}^n$ such that $K\subset K_0\subset \mathbb{C}^n\setminus\bigcup_{l\in\mathbb{N}} H_l$.
\end{enumerate}
By Kallin's lemma (see E.\ Kallin \cite{Kallin1965} or E.L.\ Stout \cite[p.\ 62]{Stout2007PM}) and the Oka-Weil theorem (see \cite[Theorem 1.5.1]{Stout2007PM}), if $C$ and $T$ are disjoint compact convex sets in $\mathbb{C}^n$ and $A\subset C$ is a compact polynomially convex set, then the union $A\cup T$ is polynomially convex as well; we refer to  \cite{Stout2007PM} for a monograph on polynomial convexity. This easily implies in view of conditions (a), (b), and (c) that:
\begin{enumerate}[(d)]
\item[(d)] The set $K_j\cup L_j$ is polynomially convex for all $j\in\mathbb{N}\cup\{0\}$. More generally, the compact set $K_j\cup \bigcup_{i=j}^k L_k$ is polynomially convex for all $k\ge j$, $j\in\mathbb{N}\cup\{0\}$.
\end{enumerate}
In this situation, the Anders\'en-Lempert theory (see E.\ Anders\'en and L.\ Lempert \cite{AndersenLempert1992IM} and F.\ Forstneri\v c and J.-P.\ Rosay \cite{ForstnericRosay1993IM}) enables one to construct a sequence $\{\Phi_j\}_{j\in\mathbb{N}}$ of holomorphic automorphisms of $\mathbb{C}^n$ such that if
\[
	\Psi_j=\Phi_j\circ\cdots\circ\Phi_1,\quad j\in\mathbb{N},
\]
then $|\Phi_j(z)-z|<\epsilon_j$ for all $z\in K_j$ and 
\[
	\Psi_j(Z)\cap\bigcup_{i=1}^j L_j=\varnothing\quad \text{for all }j\in\mathbb{N}.
\] 
Here, each $\epsilon_j$ is a positive number which may be chosen arbitrarily small at each step in the inductive construction. In fact, setting $\Psi_0=\Phi_0={\rm Id}_{\mathbb{C}^n}$ and assuming that we have such automorphisms $\Phi_1,\ldots,\Phi_{j-1}$ for some $j\in\mathbb{N}$, up to a suitable affine complex change of coordinates we can even choose $\Phi_j$ to be a shear map of the form
\[
	\Phi_j(\zeta,w)=(\zeta,e^{\phi_j(\zeta)}w),\quad \zeta\in\mathbb{C},\; w\in\mathbb{C}^{n-1},
\]
with a holomorphic function $\phi_j\colon\mathbb{C}\to\mathbb{C}$. Choosing the number $\epsilon_j>0$ sufficiently small at each step, we have that $\lim_{j\to\infty}\Psi_j=\Psi$ exists uniformly on compact subsets in
$\Omega:=\bigcup_{j\in\mathbb{N}}\Psi_j^{-1}(K_j)$,
and $\Psi\colon\Omega\to\mathbb{B}_n$ is a biholomorphic map such that 
\[
	\Psi(Z\cap\Omega)\cap L=\varnothing.
\] 
This implies that $X=Z\cap\Omega$ is a pseudoconvex Runge domain in $Z$, $K\subset X$, and the closed complex hypersurface $V=\Psi(X)$ of $\mathbb{B}_n$ is complete (since it is disjoint from the labyrinth $L$ in \eqref{eq:Lab}; see condition (ii) in the description of $L$), thereby proving Theorem \ref{th:AGL}.

Despite the fact that the method developed in \cite{AlarconGlobevnikLopez2019Crelle} works fine if one uses the labyrinths constructed in \cite{Globevnik2015AM}, the authors introduced in \cite{AlarconGlobevnikLopez2019Crelle} new labyrinths of compact sets in the ball with properties analogous to those in (i)-(iv) and (a)-(d) above. The main reason for building these new labyrinths is that they are constructed in a much simpler way; only basic trigonometry is required.
%
%
\begin{definition}\label{def:L}
Given an integer $n\ge 2$ and a pair of positive numbers $0<r<R$, we say that a compact set $\mathcal{L}$ in the spherical shell $R\mathbb{B}_n\setminus r\overline{\mathbb{B}}_n=\{z\in\mathbb{C}^n\colon r<|z|<R\}$ is a tangent labyrinth if the following conditions are satisfied:
\begin{itemize}
\item $\mathcal{L}$ has finitely many connected components $T_1,\ldots,T_k$ $(k\in\mathbb{N})$.
\smallskip
\item Each component $T_j$ of $\mathcal{L}$ is a closed round ball in a real affine hyperplane in $\mathbb{C}^n$ which is orthogonal to the position vector of the center $x_j$ of the ball $T_j$.
\smallskip
\item If $|x_i|=|x_j|$ for some $i,j\in\{1,\ldots,k\}$, then the radii of the balls $T_i$ and $T_j$ are equal, while if $|x_i|<|x_j|$ then $T_i\subset |x_j|\mathbb{B}_n$.
\end{itemize}
\end{definition}

\begin{figure}[ht]
\includegraphics[scale=.25]{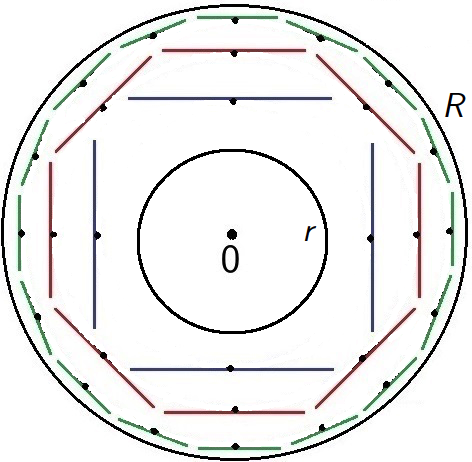}
\caption{A tangent labyrinth in the spherical shell $R\mathbb{B}_n\setminus r\overline{\mathbb{B}}_n$. Each piece of the labyrinth is a closed round ball in a real affine hyperplane. The pieces of the labyrinth are organized in layers according to the norm of the center of each piece.}
\label{fig:tidy} 
\end{figure}

These new labyrinths are much more symmetric than Globevnik's ones in \cite{Globevnik2015AM}, and its pieces are all round balls. See Figure \ref{fig:tidy}.
%
%
The main two properties of tangent labyrinths are the following:
\begin{enumerate}[(II)]
\item[(I)] If $\mathcal{L}$ is a tangent labyrinth in the spherical shell $R\mathbb{B}_n\setminus r\overline{\mathbb{B}}_n$ $(0<r<R,\; n\ge 2)$, then the compact set $r\overline{\mathbb{B}}_n\cup\mathcal{L}\subset\mathbb{C}^n$ is polynomially convex.
\smallskip
\item[(II)] For any given numbers $0<r<R$, $\delta>0$ (big), and $\tau>0$ (small) there exists a tangent labyrinth $\mathcal{L}$ in $R\mathbb{B}_n\setminus r\overline{\mathbb{B}}_n$ $(n\ge 2)$ enjoying the following conditions:
\begin{itemize}
\item If $\gamma\colon[0,1]\to\mathbb{B}_n$ is a path crossing the shell (i.e., $|\gamma(0)|\le r$ and $|\gamma(1)|\ge R$) and $\gamma$ is disjoint from the labyrinth (i.e., $\gamma([0,1])\cap\mathcal{L}=\varnothing$), then the Euclidean length of $\gamma$ is greater than $\delta$.
\smallskip
\item Each component $T$ of $\mathcal{L}$ has Euclidean diameter smaller than $\tau$; i.e. the radius of the Euclidean ball $T$ is smaller than $\tau/2$.
\end{itemize}
\end{enumerate}
Condition (I) can be seen as a consequence of the Kallin lemma and the Oka-Weil theorem, as above; this was observed by A.\ Alarc\'on and F.\ Forstneri\v c in \cite{AlarconForstneric2017PAMS}, see also \cite[Remark 2.4]{Alarcon2022JDG}. On the other hand, condition (II) was granted in \cite[Lemma 2.4]{AlarconGlobevnikLopez2019Crelle}; see also \cite[Lemma 2.3]{Alarcon2022JDG}. 

Note that the tangent labyrinths $\mathcal{L}$ in Definition \ref{def:L} are compact and have finitely many components, while those $L$ in \eqref{eq:Lab} are closed in the ball and have infinitely many components. It is however obvious in view of condition (II) that for any increasing sequence $0<r_1<R_1<r_2<R_2<\cdots$ with $\lim_{j\to\infty}r_j=1$ there is a sequence of tangent labyrinths $\mathcal{L}_j$ in $R_j\mathbb{B}_n\setminus r_j\overline{\mathbb{B}}_n$, $j\in\mathbb{N}$, with the property that every divergent path $\gamma\colon[0,1)\to\mathbb{B}_n$ meeting at most finitely many components of the infinite labyrinth 
\begin{equation}\label{eq:infiniteL}
	\mathcal{L}=\bigcup_{j\in\mathbb{N}}\mathcal{L}_j\subset\mathbb{B}_n
\end{equation} 
has infinite length; cf.\ property (ii) of the Globevnik labyrinth $L$ in \eqref{eq:Lab}.

The new construction method developed in \cite{AlarconGlobevnikLopez2019Crelle}, based on the use of tangent labyrinths and holomorphic automorphisms of $\mathbb{C}^n$, has given rise to some further developments. For instance, A.\ Alarc\'on and J.\ Globevnik extended Corollary \ref{co:AGL} to complex curves of arbitrary topology, possibly infinite, including in addition a hitting condition.
%
%
\begin{theorem}\label{th:AG-C2}
{\bf (Alarc\'on-Globevnik \cite{AlarconGlobevnik2017C2})} Given a closed discrete subset $\Lambda$ of $\mathbb{B}_2$ and an open orientable smooth surface $M$, there is a complex structure $J$ on $M$ such that the open Riemann surface $(M,J)$ admits a complete proper holomorphic embedding $(M,J)\hookrightarrow \mathbb{B}_2$ whose image contains $\Lambda$.
\end{theorem}

Choosing the set $\Lambda$ such that the limit set $\overline{\Lambda}\setminus\Lambda$ equals the boundary sphere $b\mathbb{B}_n=\{z\in\mathbb{C}^n\colon |z|=1\}$, one obtains the following immediate corollary.
%
%
\begin{corollary}\label{co:AG-C2}
{\bf (Alarc\'on-Globevnik \cite{AlarconGlobevnik2017C2})} $\mathbb{B}_2$ contains complete properly embedded complex curves with any topology whose limit set is the whole sphere $b\mathbb{B}_2$.
\end{corollary}

In the case when $M$ is simply-connected, and hence the Riemann surface $(M,J)$ is biholomorphic to the disc $\mathbb{D}$, Theorem \ref{th:AG-C2} was previously obtained by J.\ Globevnik in \cite{Globevnik2016JMAA}. The proof of Theorem \ref{th:AG-C2} combines ideas from \cite{AlarconGlobevnikLopez2019Crelle} and \cite{AlarconLopez2013JGEA,ForstnericWold2009}, the latter relying on the method of exposing boundary points of bordered Riemann surfaces in $\mathbb{C}^2$ as well as on the use of Fatou-Bieberbach domains of $\mathbb{C}^2$. Theorem \ref{th:AG-C2} is somehow analogous to the main result in \cite{AlarconLopez2013JGEA} where properly embedded complex curves in $\mathbb{C}^2$ with arbitrary topology are constructed.

Also using tangent labyrinths and holomorphic automorphisms of $\mathbb{C}^n$ as in \cite{AlarconGlobevnikLopez2019Crelle}, A.\ Alarc\'on and F.\ Forstneri\v c obtained complete complex submanifolds in $\mathbb{B}_n$ with controlled topology, having no self-intersections, and hitting any given countable subset (not necessarily closed or discrete), thereby complementing Theorem \ref{th:AGL}.
%
%
\begin{theorem}\label{th:AF2017}
{\bf (Alarc\'on-Forstneri\v c \cite{AlarconForstneric2017PAMS})}
Let $\Lambda\subset \mathbb{B}_n$ $(n\ge 2)$ be a countable set.
If $Z$ is a smooth closed complex submanifold of $\mathbb{C}^n$ and $K\neq\varnothing$ is a connected compact subset of $Z\cap\mathbb{B}_n$, then there are a pseudoconvex Runge domain $X\subset Z$ such that $K\subset X$ and a complete holomorphic injective immersion $\varphi\colon X\hookrightarrow\mathbb{B}_n$ with $\Lambda\subset\varphi(X)$. 
In particular, $\varphi(X)$ can be made dense in $\mathbb{B}_n$.
\end{theorem}

In the case of complex curves, this has the following immediate corollary; compare with Theorem \ref{th:AG-C2} and statement 3 in Theorem \ref{th:Jordan}.
%
%
\begin{corollary}\label{co:AF2017}
{\bf (Alarc\'on-Forstneri\v c \cite{AlarconForstneric2017PAMS})}
Given a countable subset $\Lambda$ of $\mathbb{B}_n$ $(n\ge 2)$ and an open orientable smooth surface $M$ of finite topology, there is a complex structure $J$ on $M$ such that the open Riemann surface $(M,J)$ admits a complete holomorphic injective immersion $\varphi\colon (M,J)\hookrightarrow \mathbb{B}_n$ whose image contains $\Lambda$.

In particular, $\varphi(M)$ can be made dense in $\mathbb{B}_n$.
\end{corollary}

The first result which provided control of the topology of all leaves in a holomorphic foliation of the ball by complete closed complex submanifolds was subsequently given by A.\ Alarc\'on and F.\ Forstneri\v c.
%
%
\begin{theorem}\label{th:AF2020}
{\bf (Alarc\'on-Forstneri\v c \cite{AlarconForstneric2020MZ})}
There exists a nonsingular holomorphic foliation of the ball $\mathbb{B}_n$ $(n\ge 2)$ all of whose leaves are complete properly embedded complex discs in $\mathbb{B}_n$.
\end{theorem}

The topology of the leaves in the foliations given in Theorem \ref{th:AF2020} is the simplest possible one: discs. It is not known whether a comparable result holds for leaves with prescribed but more complicated topology. The proof follows a similar approach as that in \cite{AlarconGlobevnikLopez2019Crelle}, but with a more careful use of the Anders\'en-Lempert theory. The idea is, using holomorphic automorphisms, to successively twist a holomorphic foliation of $\mathbb{C}^n$ by complex lines in order to make bigger and bigger parts of the foliation avoid more and more pieces of an infinite labyrinth in $\mathbb{B}_n$, as those in  \eqref{eq:Lab} or \eqref{eq:infiniteL}. Each component of the intersection of the ball with a properly embedded complex line is Runge in the complex line, and hence simply connected. Since it is in addition bounded, it is a properly embedded complex disc in the ball. More precisely, the construction provides a sequence of foliations $\{\mathcal{F}_j\}_{j\in\mathbb{N}}$ of $\mathbb{B}_n$ by properly embedded complex discs, converging to a limit foliation $\mathcal{F}$ whose all leaves are discs and each one meets at most finitely many components of the labyrinth, so it is complete.

None of the methods discussed in this section enables one to control the complex structure of the examples, except of course if the curve is simply connected. In particular, the following remains an open question; see \cite[Problem 1.5]{AlarconGlobevnik2017C2}.
%
%
\begin{problem}\label{pr:compact}
Let $\overline M=M\cup bM$ be a compact bordered Riemann surface admitting a smooth embedding $\overline M\hookrightarrow\mathbb{C}^2$ which is holomorphic on $M$. Does there exist a complete holomorphic embedding $M\hookrightarrow\mathbb{C}^2$ with bounded image?
\end{problem}

F.\ Forstneri\v c and E.F.\ Wold proved in \cite{ForstnericWold2009} that every bordered Riemann surface $M$ as in Problem \ref{pr:compact} admits a proper holomorphic embedding $M\hookrightarrow\mathbb{C}^2$. Obviously, such an embedding is complete but unbounded. Here is a less ambitious but still appealing question in this direction. 
%
%
\begin{problem}
Let $\overline M$ be as in Problem  \ref{pr:compact}. Does there exist a complete non-proper holomorphic embedding $M\hookrightarrow\mathbb{C}^2$?
\end{problem}
%
%
%

%
%
\subsection{Complete complex hypersurfaces in the ball come in foliations}

In principle, it was not clear whether the construction of Globevnik in \cite{Globevnik2015AM} gives, say, a rare type or a very special type of solution to the embedded Yang problem for hypersurfaces in the ball. Note that his examples of complete closed complex hypersurfaces in $\mathbb{B}_n$ are obtained as leaves of a holomorphic foliation of $\mathbb{B}_n$ by complex hypersurfaces of the same sort. Thus, the question whether every such hypersurface can be embedded into a foliation of the ball by hypersurfaces with the same properties naturally appeared. Recall that, at that time, there were already available in the literature several constructions of complete properly embedded complex submanifolds in balls, besides the implicit one, so one could be inclined to think that the answer to this question is negative. However, the answer is positive, as it follows from the following more general result by A.\ Alarc\'on dealing with submanifolds of arbitrary codimension.
%
%
\begin{theorem}\label{th:JDG}
{\bf (Alarc\'on \cite{Alarcon2022JDG})}
Let $n$ and $q$ be integers with $1\le q<n$, let $V$ be a smooth closed complex submanifold of pure codimension $q$ in $\mathbb{B}_n$, possibly disconnected, and assume that $V$ is contained in a fibre of a holomorphic submersion $\mathbb{B}_n\to\mathbb{C}^q$. Then there is a holomorphic submersion $f\colon\mathbb{B}_n\to\mathbb{C}^q$ such that:
\begin{itemize}
\item $V$ is a union of components of $f^{-1}(0)$.
\smallskip
\item The fibre $f^{-1}(c)\subset\mathbb{B}_n$ is a smooth complete closed complex submanifold of pure codimension $q$ for every $c\in f(\mathbb{B}_n)\setminus\{0\}$.
\smallskip
\item $f^{-1}(0)\setminus V$ is either empty or a smooth complete closed complex submanifold of pure codimension $q$ in $\mathbb{B}_n$.
\smallskip
\item If $V$ is a fibre of a holomorphic submersion $\mathbb{B}_n\to\mathbb{C}^q$, then $f$ can be chosen with $f^{-1}(0)=V$.
\end{itemize}
\end{theorem}

The family of components of the fibres $f^{-1}(c)$ $(c\in\mathbb{C}^q)$ of the holomorphic submersion $f$ given by the theorem form a nonsingular holomorphic submersion foliation of $\mathbb{B}_n$ by smooth connected closed complex submanifolds of codimension $q$ all which, except perhaps those contained in the initial submanifold $V\subset f^{-1}(0)$, are complete. If the given $V$ is complete, then all the leaves in the foliation are complete. This provided the first known example of a foliation of the ball $\mathbb{B}_n$ by complete connected closed complex submanifolds of codimension $q>1$, as well as the first example of a nonsingular holomorphic foliation of $\mathbb{B}_n$ by such submanifolds for any dimension $n\ge 2$ and codimension $1\le q<n$. 

Somewhat surprisingly, the hypersurface $V$ in Theorem \ref{th:JDG} need not be complete, so the following holds.
%
%
\begin{corollary}\label{co:JDG}
{\bf (Alarc\'on \cite{Alarcon2022JDG,Alarcon2022IUMJ})}
Let $n$ and $q$ be a pair of integers with $1\le q<n$. The following assertions hold: 
\begin{enumerate}[1.]
\item There exists a nonsingular holomorphic submersion foliation of $\mathbb{B}_n$ by smooth connected complete closed complex submanifolds of codimension $q$.
\smallskip
\item For any integer $k\ge 1$, there exists a nonsingular holomorphic submersion foliation of $\mathbb{B}_n$ by smooth connected closed complex submanifolds of codimension $q$ all which, except precisely $k$ among them, are complete. 
\smallskip
\item There exists a nonsingular holomorphic submersion foliation of $\mathbb{B}_n$ by smooth connected closed complex submanifolds of codimension $q$ all which, except precisely countably infinitely many, are complete.
\end{enumerate}
\end{corollary}

The following result concerning the existence of foliations of the ball having both many complete leaves and many incomplete leaves was subsequently obtained by a generalization of the construction method in  \cite{Alarcon2022JDG}.
%
%
\begin{corollary}\label{co:Indiana}
{\bf (Alarc\'on \cite{Alarcon2022IUMJ})}
Given integers $n$ and $q$ with $1\le q<n$, there is a nonsingular holomorphic submersion foliation $\mathcal{F}$ of the ball $\mathbb{B}_n$ by smooth connected closed complex submanifolds of codimension $q$ such that both the union of the complete leaves of $\mathcal{F}$ and the union of the incomplete leaves of $\mathcal{F}$ are dense subsets of $\mathbb{B}_n$. In particular, every leaf of $\mathcal{F}$ is both a limit leaf of the family of complete leaves of $\mathcal{F}$ and a limit leaf of the family of incomplete leaves of $\mathcal{F}$.
\end{corollary}

The assumption in Theorem \ref{th:JDG} that the given hypersurface $V$ is contained in a fibre of a holomorphic submersion is necessary and cannot be relaxed. In the case of $q=1$, this condition is always satisfied (see F.\ Forstneri\v c \cite{Forstneric2018PAMS}), so the theorem applies to every smooth closed complex hypersurface in $\mathbb{B}_n$. More precisely, it is known that every closed complex hypersurface $V$ (possibly with singularities) in $\mathbb{B}_n$ is defined by a holomorphic function on $\mathbb{B}_n$ which is noncritical out of the singular set of $V$. (This actually holds true with the ball replaced by any Stein manifold $X$ with the vanishing second cohomology group $H^2(X;Z)=0$  \cite{Forstneric2018PAMS}; which extends the classical result by J.-P.\ Serre from 1953 \cite{Serre1953} that every divisor in such a Stein manifold $X$ is a principal divisor.) The following more precise version of Theorem \ref{th:JDG} in the case of hypersurfaces, possibly with singularities, was also obtained in \cite{Alarcon2022JDG}.
%
%
\begin{theorem}\label{th:JDG-H}
{\bf (Alarc\'on \cite{Alarcon2022JDG})}
For any closed complex hypersurface $V$ (possibly with singularities) in $\mathbb{B}_n$ $(n\ge 2)$ there exists a holomorphic function $f$ on $\mathbb{B}_n$ satisfying the following conditions:
\begin{itemize}
\item $V=f^{-1}(0)$.
\smallskip
\item The critical locus of $f$ coincides with the singular set of $V$.
\smallskip
\item The level set $f^{-1}(c)\subset \mathbb{B}_n$ is a smooth complete closed complex hypersurface for every $c\in f(\mathbb{B}_n)\setminus\{0\}$.
\end{itemize}
\end{theorem}

In this case, the family of components of the level sets $f^{-1}(c)$ $(c\in\mathbb{C})$ of the holomorphic function $f$ furnished by the theorem form a holomorphic foliation of $\mathbb{B}_n$ by connected closed complex hypersurfaces all which, except perhaps those contained in the given hypersurface $V= f^{-1}(0)$, are smooth and complete. If the given $V$ is smooth, then $f$ is noncritical and hence the foliation nonsingular; while if $V$ is complete, then all the leaves in the foliation are complete. This established a sort of converse to Theorem \ref{th:PikoMA}, namely, it proves that, in fact, every complete closed complex hypersurface in $\mathbb{B}_n$ is defined by a holomorphic function on $\mathbb{B}_n$ whose level sets are all complete (i.e., as those constructed by Globevnik in \cite{Globevnik2015AM}), and hence it can be embedded as a leaf in a holomorphic foliation of the ball by complete closed complex hypersurfaces. Moreover, the defining function can be chosen noncritical, and hence the foliation nonsingular, whenever the given hypersurface is smooth.

The method used in the proof of Theorem \ref{th:JDG} broadly follows the one developed by Globevnik for proving Theorem \ref{th:PikoMA} in \cite{Globevnik2015AM}, but it presents some major differences and novelties. The holomorphic submersion $f\colon\mathbb{B}_n\to\mathbb{C}^q$ satisfying the conclusion of the theorem is obtained as the limit of a sequence of holomorphic submersions $f_j\colon\mathbb{B}_n\to\mathbb{C}^q$ $(j\in\mathbb{N})$, and it is granted that each fibre of $f$, except perhaps the one over $0\in\mathbb{C}^q$, meets at most finitely many components of a suitable infinite labyrinth of compact sets in $\mathbb{B}_n$, thereby making sure that every such fibre is complete. A main novelty is that the proof requires to construct the labyrinth and the limit holomorphic submersion $f$ at the same time in an inductive procedure; that is, the labyrinth cannot be fixed beforehand, as in the proofs in \cite{Globevnik2015AM,AlarconGlobevnikLopez2019Crelle} and in the aforementioned subsequent papers. Moreover, the labyrinth depends on the given submanifold $V$, and the components of the labyrinth that meet $V$ are treated in a completely different way than those that do not; in particular, the former components must be chosen with small diameter. Another novelty of the proof is that it relies on a rather sophisticated Oka-Weil-Cartan type approximation with interpolation theorem for holomorphic submersions due to F.\ Forstneri\v c (see \cite{Forstneric2003AM,Forstneric2018PAMS} or \cite[\textsection 9.12-9.16]{Forstneric2017E}), instead of the more basic Runge-Weierstrass theorem. We next spell out some of the details.

Let us briefly outline the proof of Theorem \ref{th:JDG} in the slightly simpler case when there is a holomorphic submersion $f_0\colon\mathbb{B}_n\to\mathbb{C}^q$ such that $V=f_0^{-1}(0)$; that is always the case if $q=1$. Choose any increasing sequence of positive numbers $0<r_1<R_1<r_2<R_2<\cdots\to 1$, fix $\epsilon_0>0$, choose $0<R_0<r_1$, and call $L_0=\varnothing$. The proof consists of inductively constructing a sequence $S_j=\{f_j,\epsilon_j,\mathcal{L}_j\}_{j\in\mathbb{N}}$, where
\begin{itemize}
\item $f_j\colon\mathbb{B}_n\to\mathbb{C}^q$ is a holomorphic submersion,
\smallskip
\item $\epsilon_j>0$ is a positive number, and
\smallskip
\item $\mathcal{L}_j$ is a tangent labyrinth in the spherical shell $R_j\mathbb{B}_n\setminus r_j\overline{\mathbb{B}}_n$ (see Definition \ref{def:L}),
\end{itemize}
such that the following conditions are satisfied for each $j\in\mathbb{N}$:
\begin{enumerate}[(Aa)]
\item[(1$_j$)] $|f_j(z)-f_{j-1}(z)|<\epsilon_j$ for all $z\in r_j\overline{\mathbb{B}}_n$.
\smallskip 
\item[(2$_j$)] $f_j^{-1}(0)=V$.
\smallskip 
\item[(3$_j$)] $0<2\epsilon_j<\epsilon_{j-1}$ and if $f\colon\mathbb{B}_n\to\mathbb{C}^q$ is a holomorphic map such that $|f(z)-f_{j-1}(z)|<2\epsilon_j$ for all $z\in r_j\overline{\mathbb{B}}_n$, then $f$ is submersive everywhere on $R_{j-1}\overline{\mathbb{B}}_n$.
\smallskip 
\item[(4$_j$)] If $\gamma\colon[0,1]\to\mathbb{B}_n$ is a path satisfying that $|\gamma(0)|\le r_j$, $|\gamma(1)|\ge R_j$, and $\gamma([0,1])\cap\mathcal{L}_j=\varnothing$, then ${\rm length}(\gamma)>1$.
\smallskip 
\item[(5$_j$)] If $z\in\mathcal{L}_j$, then either $|f_j(z)|>j$ or $|f_j(z)|<1/j$.
\end{enumerate}
Carrying out this process in the right way, one obtains a holomorphic submersion 
\[
	f=\lim_{j\to\infty}f_j\colon\mathbb{B}_n\to\mathbb{C}^q,\quad j\in\mathbb{N},
\]
as well as an infinite labyrinth of compact sets
\[
	\mathcal{L}=\bigcup_{j\in\mathbb{N}}\mathcal{L}_j\subset\mathbb{B}_n,
\]
satisfying the following conditions:
\begin{enumerate}[(B)]
\item[(A)] $f^{-1}(0)=V$.
\smallskip
\item[(B)] Every divergent path $[0,1)\to\mathbb{B}_n$ hitting at most finitely many components of $\mathcal{L}$ has infinite length.
\smallskip
\item[(C)] $f^{-1}(c)$ meets at most finitely many components of $\mathcal{L}$ for each $c\in\mathbb{C}^q\setminus\{0\}$.
\end{enumerate}
Note that conditions (B) and (C) are granted by properties (4$_j$) and (5$_j$), respectively. 
It is clear that (B) and (C) ensure that $f^{-1}(c)$ is a smooth complete closed complex submanifold of pure codimension $q$ in $\mathbb{B}_n$ for all $c\in\mathbb{C}^q\setminus\{0\}$, and hence, in view of condition (A), $f$ satisfies the conclusion of the theorem. 

For the inductive step, fix $j\in\mathbb{N}$ and assume that we already have a suitable triple $S_{j-1}=\{f_{j-1},\epsilon_{j-1},\mathcal{L}_{j-1}\}$. First of all, choose a number $\epsilon_j>0$ so small that (3$_j$) holds; that is always possible by the Cauchy estimates. Also choose, as we may by (2$_{j-1}$), a number $\tau>0$ so small that
\begin{equation}\label{eq:tau}
	|f_{j-1}(z)|<1/j\quad \text{for all $z\in R_j\overline{\mathbb{B}}_n$ with ${\rm dist}(z,V)<\tau$}.
\end{equation}
Next, take a tangent labyrinth $\mathcal{L}_j$ in $R_j\mathbb{B}_n\setminus r_j\overline{\mathbb{B}}_n$ satisfying (4$_j$) and such that
\begin{equation}\label{eq:tau-1}
	{\rm diam}(T)<\tau\quad \text{for each component $T$ of $\mathcal{L}_j$}.
\end{equation}
(See property (II) below Definition \ref{def:L}.) We then split the labyrinth $\mathcal{L}_j$ into two disjoint parts, $\mathcal{L}_j=\Lambda_V\cup\Lambda_0$, where $\Lambda_V$ is the union of those components of $\mathcal{L}_j$ which intersect $V$ and $\Lambda_0$ is the union of those which are disjoint from $V$. Both $\Lambda_V$ and $\Lambda_0$ consist of finitely many components and $\Lambda_V\cap\Lambda_0=\varnothing$. So, we have in view of \eqref{eq:tau} and \eqref{eq:tau-1} that 
\begin{equation}\label{eq:LambdaV}
	|f_{j-1}(z)|<1/j\quad \text{for all $z\in \Lambda_V$}.
\end{equation}
In this situation, a standard application of the aforementioned Oka-Weil-Cartan-Forstneri\v c theorem for holomorphic submersions (see \cite[Theorem 2.1]{Alarcon2022JDG} for the precise statement that we need here), taking into account that $r_j\overline{\mathbb{B}}_n\cup \mathcal{L}_j$ is a polynomially convex compact set (see property (I) below Definition \ref{def:L}), furnishes us with a holomorphic submersion $f_j\colon \mathbb{B}_n\to\mathbb{C}^q$ satisfying the following conditions:
\begin{itemize}
\item $|f_j(z)-f_{j-1}(z)|<\delta$ for all $z\in r_j\overline{\mathbb{B}}_n\cup\Lambda_V$ for any given $\delta>0$.
\smallskip
\item $f_j^{-1}(0)=V$.
\smallskip
\item $|f_j(z)|>j$ for all $z\in \Lambda_0$.
\end{itemize}
It is clear that these properties, together with \eqref{eq:LambdaV} and the fact that $\mathcal{L}_j=\Lambda_V\cup\Lambda_0$, guarantee that $f_j$ satisfies conditions (1$_j$), (2$_j$), and (5$_j$) provided that the number $\delta>0$ is chosen sufficiently small. This closes the induction and completes the sketch of the proof of Theorem \ref{th:JDG}.

Complete complex hypersurfaces in a Stein manifold equipped with a Riemannian metric come in foliations as well, as the following result points out.
%
%
\begin{theorem}\label{th:JDG-Stein}
{\bf (Alarc\'on \cite{Alarcon2022JDG})}
Let $X$ be a Stein manifold of dimension $\ge 2$ equipped with a Riemannian metric, and assume that $V$ is a smooth closed complex hypersurface in $X$ such that the normal bundle $N_{V/X}$ of $V$ in $X$ is trivial. Then there exists a holomorphic function $f$ on $X$ satisfying the following conditions:
\begin{itemize}
\item $V\subset f^{-1}(0)$.
\smallskip
\item The level set $f^{-1}(c)$ is complete for every $c\in f(X)\setminus\{0\}$.
\smallskip
\item $f^{-1}(0)\setminus V$ is either empty or a complete closed complex hypersurface (possibly with singularities) in $X$.
\smallskip
\item $f$ in noncritical everywhere on $V$.
\end{itemize}
\end{theorem}

Completeness in this statement refers to the metric induced in the hypersurface by that in $X$. Note that the assumption on the normal bundle $N_{V/X}$ is necessary by the F.\ Docquier and H.\ Grauert tubular neighborhood theorem \cite{DocquierGrauert1960MA}. The connected components of the level sets $f^{-1}(c)\subset X$ $(c\in\mathbb{C})$ of the holomorphic function $f$ provided by Theorem \ref{th:JDG-Stein} form a (possibly singular) holomorphic foliation of the Stein manifold $X$ by connected closed complex hypersurfaces all which, except perhaps those contained in $V$, are complete. If $V$ is complete, then all the leaves are complete. Moreover, most level sets of $f$ are smooth by Sard's theorem, and hence we obtain the following extension of Corollary \ref{co:PikoMA}. (Recall that pseudoconvex domains are Stein manifolds, but not every Stein manifold is biholomorphic to a pseudoconvex domain in a complex Euclidean space.)
%
%
\begin{corollary}\label{co:JDG-Stein}
{\bf (Alarc\'on \cite{Alarcon2022JDG})}
Every Stein manifold $X$ of dimension $\ge 2$ equipped with a Riemannian metric admits a (possibly singular) holomorphic foliation by complete connected closed complex hypersurfaces (possibly with singularities).

In particular, $X$ admits a complete properly embedded complex hypersurface.
\end{corollary}

The proof of Theorem \ref{th:JDG-Stein} combines the construction method in \cite{Alarcon2022JDG} with a slight refinement of the trick used by Globevnik in \cite{Globevnik2016MA} to show that every pseudoconvex domain (other than the ball) in $\mathbb{C}^n$ admits a holomorphic foliation by complete closed complex hypersurfaces (see Theorem \ref{th:PikoMA}).

Although every Stein manifold admits a holomorphic function without critical points, as was shown by F.\ Forstneri\v c in \cite{Forstneric2003AM}, the method of proof in \cite{Alarcon2022JDG} does not allow to guarantee that the function $f$ furnished by Theorem \ref{th:JDG-Stein} be noncritical. Therefore, the following remains an open problem (see \cite[Conjecture 5.2]{Alarcon2022JDG}).
%
%
\begin{problem}\label{pr:Stein}
Does every Stein manifold $X$ of dimension $n\ge 2$ equipped with a Riemannian metric admit a nonsingular holomorphic foliation by smooth complete closed complex hypersurfaces? 
\end{problem}

%
%
\subsection{Labyrinths in pseudoconvex domains}

Recently, S. Charpentier and {\L}. Kosi\'{n}ski \cite{CharpentierKosinski2020}  constructed labyrinths of compact sets in any given pseudoconvex Runge domain $D\subset\mathbb{C}^n$ with properties analogous to those of the tangent labyrinths in balls (see Definition \ref{def:L} and properties (I) and (II) below it). 
%
%
\begin{theorem}\label{th:CK}
{\bf (Charpentier-Kosi\'{n}ski \cite{CharpentierKosinski2020})}
Let $\Omega\subset\mathbb{C}^n$ $(n\ge 2)$ be a smooth strongly pseudoconvex domain and let $K\subset \Omega$ be a $\mathcal{O}(\Omega)$-convex compact set. Then, for any $\delta>0$ there is a compact set $\Gamma\subset \Omega\setminus K$ satisfying the following conditions:
\begin{itemize}
\item $\Gamma$ has finitely many components all of which are holomorphically contractible.
\smallskip
\item $K\cup \Gamma$ is $\mathcal{O}(\Omega)$-convex.
\smallskip
\item If $\gamma\colon[0,1)\to\Omega$ is a divergent path with $\gamma(0)\in K$ and $\gamma([0,1))\cap \Gamma=\varnothing$, then ${\rm length}(\gamma)>\delta$.
\smallskip
\item If $\Omega$ is a Runge domain, then the components of $\Gamma$ can be chosen to be images of Euclidean balls of real dimension $2n-1$ under $\mathbb{R}$-affine isomorphisms and automorphisms of $\mathbb{C}^n$.
\end{itemize}
\end{theorem}

Since every pseudoconvex domain $D\subset\mathbb{C}^n$   admits a normal exhaustion $\Omega_1\Subset \Omega_2\Subset\cdots\subset\bigcup_{j\in\mathbb{N}}\Omega_j=D$  by smooth $\mathcal{O}(D)$-convex strongly pseudoconvex domains, which can be taken to be Runge if $D$ is Runge (see \cite[\textsection 2.3]{Forstneric2017E}), an inductive application of Theorem \ref{th:CK} provides holomorphically contractible compact sets $\Gamma_j\subset \Omega_{j+1}\setminus \overline{\Omega}_j$, $j\in\mathbb{N}$, satisfying the following conditions (see \cite[Theorem 1.1]{CharpentierKosinski2020}):
\begin{itemize}
\item $\Gamma_j$ has finitely many components and $\overline{\Omega}_j\cup \Gamma_j$ is $\mathcal{O}(D)$-convex for all $j\in\mathbb{B}_n$.
\smallskip
\item If $\gamma\colon[0,1)\to D$ is a divergent path meeting at most finitely many components of $\Gamma=\bigcup_{j\in\mathbb{N}}\Gamma_j$, then $\gamma$ has infinite length.
\smallskip
\item If $\Omega$ is a Runge domain, then the components of $\Gamma$ can be chosen to be images of Euclidean balls of real dimension $2n-1$ under $\mathbb{R}$-affine isomorphisms and automorphisms of $\mathbb{C}^n$.
\end{itemize}

This can be used to extend from the ball $\mathbb{B}_n$ to any pseudoconvex Runge domain $D$ in $\mathbb{C}^n$ many of the results mentioned in this survey, just by using the labyrinths in Theorem \ref{th:CK} instead of the tangent labyrinths in Definition \ref{def:L}. In particular, all results in Section \ref{sec:topology}, as well as Theorem \ref{th:JDG}, Theorem \ref{th:JDG-H}, and the first assertion in Corollary \ref{co:JDG}, hold true with the ball replaced by any pseudoconvex Runge domain (for the extension of Theorem \ref{th:AF2020} it is also required that the domain $D$ be hyperbolic, meaning that it does not contain any biholomorphic copy of $\mathbb{C}$). Note, in particular, that this provides an affirmative answer to Problem \ref{pr:Stein} in the special case when the Stein manifold $X$ is a pseudoconvex Runge domain in $\mathbb{C}^n$.
 Likewise, Corollary \ref{co:Indiana} and assertions 2 and 3 in Corollary \ref{co:JDG} hold true with the ball replaced by any such domain $D$ if we assume in addition that $D$ contains no smooth closed complex submanifold of $\mathbb{C}^n$; this happens for instance if $D\subset\mathbb{C}^n$ is relatively compact. 

In the subsequent paper \cite{CharpentierKosinski2021}, S. Charpentier and {\L}. Kosi\'{n}ski extended Globevnik's existence result for highly oscillating holomorphic functions in Theorem \ref{th:PikoMA}, to the existence of functions with an even wilder asymptotic behavior.
%
%
\begin{theorem}\label{th:CK-wild}
{\bf (Charpentier-Kosi\'{n}ski \cite{CharpentierKosinski2021})}
Let $D\subset\mathbb{C}^n$ $(n\ge 2)$ be a pseudoconvex domain. Then there are holomorphic functions $f$ on $D$ with the property that $f(\gamma([0,1)))$ is everywhere dense in $\mathbb{C}$ for any divergent path $\gamma\colon[0,1)\to D$ with finite length. In fact, the set of all such functions is a residual, densely lineable, spaceable subset of the space of all holomorphic functions on $D$ endowed with the locally uniform convergence topology.
\end{theorem}

Note that the functions $f$ provided by Theorem \ref{th:CK-wild} also satisfy the conclusion of Theorem \ref{th:PikoMA}. In particular, the components of the level sets $f^{-1}(c)$ $(c\in\mathbb{C})$ form a (possibly singular) holomorphic foliation of $D$ by complete connected closed complex hypersurfaces (possibly with singularities), as those in Corollary \ref{co:PikoMA}. We emphasize that the functions $f$ in Theorem \ref{th:CK-wild} are abundant; in particular, they form a dense subset of the space of all holomorphic functions on $D$. Therefore, the functions in Theorem \ref{th:PikoMA} are abundant as well.


\subsection*{Acknowledgements}
This research was partially supported by the State Research Agency (AEI) via the grant no.\ PID2020-117868GB-I00 and the ``Maria de Maeztu'' Excellence Unit IMAG, reference CEX2020-001105-M, funded by MCIN/AEI/10.13039/501100011033/; and the Junta de Andaluc\'ia grant no. P18-FR-4049; Spain. 
The author wishes to thank Franc Forstneri\v c for helpful suggestions which led to improve the exposition.




\medskip
\noindent Antonio Alarc\'{o}n

\noindent Departamento de Geometr\'{\i}a y Topolog\'{\i}a e Instituto de Matem\'aticas (IMAG), Universidad de Granada, Campus de Fuentenueva s/n, E--18071 Granada, Spain.

\noindent  e-mail: {\tt alarcon@ugr.es}

\end{document}